\documentclass[12pt]{amsart}

\usepackage[a4paper]{geometry}

\usepackage{amsfonts, amsmath, amssymb, amscd, latexsym, graphicx, pb-diagram, caption, mathdots}

\def\qed{\hfill\rlap{$\sqcup$}$\sqcap$\par}

\newtheorem{thm}{Theorem}[section]

\newtheorem{defn}[thm]{Definition}

\newtheorem{question*}{Question}

 \theoremstyle{definition}

\title{A note on the metric of Thompson's group $V$}
\author{Jos\'e Burillo}
\address{Departament de Matem\`atiques,
Universitat Polit\`ecnica de Catalunya, Jordi Girona 1-3, 08034
Barcelona, Spain} \email{pep.burillo@upc.edu}
\thanks{The first author thanks the Spanish Ministry MICINN through grant PID2021-126851NB-I00P for their support. }

\author{Marc Felipe}
\email{marc.felipe.alsina@gmail.com}

\begin{document}

\maketitle

\begin{abstract}
In this short note, a bound on the word metric for Thompson's group
$V$ given by Birget in 2004 is improved to a new bound, which agrees
with the known bounds for Thompson's groups $F$ and $T$.
\end{abstract}

\section*{Introduction }

Thompson's groups have received a large amount of research due to
their interesting properties (\cite{cfp}). A line of research on
them has been the study of their metric properties \cite{birget},
\cite{burillo}, \cite{bcs}, which give some insight on the word
metric of the groups and also provide methods of computing
distortions of some groups into others.

The estimate of the metric for Thompson's group $V$ was provided by
Birget in \cite{birget} in terms of the number of carets of each
tree in the minimal reduced diagram. If $N$ is this number of
carets, we have the inequality
$$
\frac{N(x)}C\leq ||x||\leq C\,N(x)\log N(x)
$$
for some universal constant $C$. For the subgroups $F$ and $T$, the
distance is equivalent to the number of carets, which imply that
they are undistorted as subgroups of $V$. By taking a larger
constant if necessary (either here or for Birget's bound for $V$),
we can assume that the constant is the same, namely, that there
exists a constant $C$ such that the following inequalities are
satisfied:
\begin{itemize}
\item For every $x\in V$,
\begin{equation}\label{eq:1} \frac{N(x)}C\leq ||x||\leq C\,N(x)\log N(x),
\end{equation}
\item for every $x\in F$ and every $x\in T$,
\begin{equation}\label{eq:2}
\frac{N(x)}C\leq ||x||\leq C\,N(x),
\end{equation}
\end{itemize}
where $N(x)$ is the number of carets in each tree of the (unique) reduced diagram for $x$.

Ideally, we could compute a quantity $D(x)$, based on the element
$x\in V$ and which necessarily must involve its permutation, which
would be equivalent, up to a multiplicative constant, to the
distance $||x||$ from $x$ to the identity. But such a quantity
remains elusive. The goal of this note is to provide an improvement
on Birget's bound and start paving the way to find such a quantity
$D$.

The reader is encouraged to check \cite{cfp} for an introduction to
Thompson's groups, in particular to $V$. It is assumed here that the
reader is familiar with the representation of elements of $V$ as
tree-pair diagrams with a permutation on its leaves, as well as
their representations as right-continuous maps on the interval.

Given $x\in V$, the number of carets $N(x)$ of the reduced diagram
provides the first approximation to the distance, but the bound
obtained by Birget is crude and not satisfactory, mainly because it
does not take the permutation into consideration. Here we will start
with a first quantity geared towards finding a better estimate of
the metric in $V$.

\section{Clusters and the new bound}

We start with the definition of the main objects we use to estimate
the metric.

\begin{defn} Let $\sigma$ be a permutation in $\mathcal S_n$. We
define a \emph{cluster} of $\sigma$ as a maximal set of consecutive
integers $\{i,i+1,\ldots,i+k\}$ inside $\{1,2,\ldots,n\}$ such that
$$
\sigma(i+1)=\sigma(i)+1,\quad\sigma(i+2),=\sigma(i)+2
\quad\ldots,\quad\sigma(i+k)=\sigma(i)+k
$$
i.e., a maximal set of consecutive integers whose images by $\sigma$
are also consecutive.
\end{defn}

The whole idea of this paper is that we can use the trees in the
diagram to organize a cluster in a particular shape of tree so that
it cancels and the diagram (and hence the bound) gets diminished.

\begin{figure}[t]
\centerline{\includegraphics{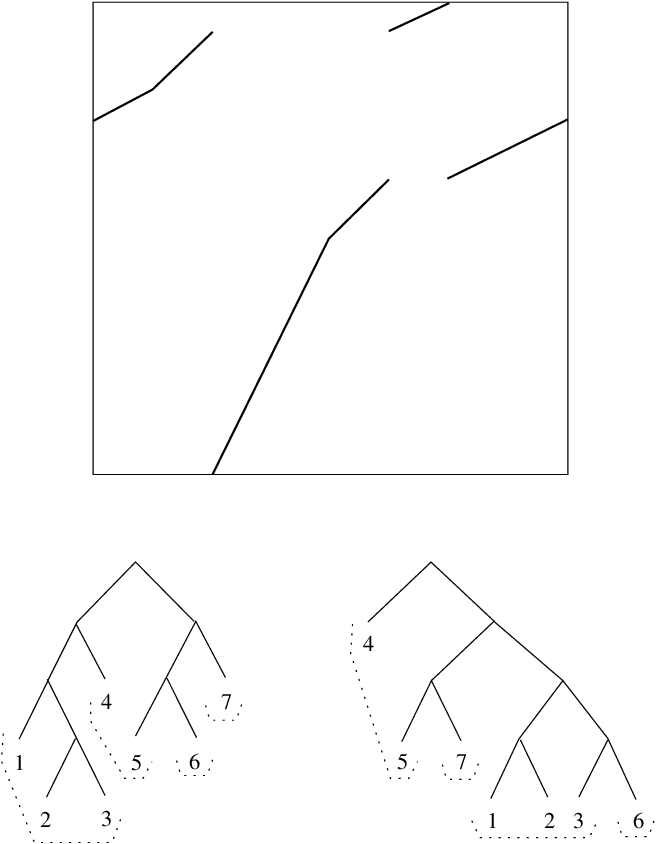}} \caption{An example of an
element of $V$ with seven leaves (six carets) but only four
clusters. The clusters are marked with dotted lines on the diagram,
and observe that they coincide with the connected components of the
graph of the corresponding map.} \label{element}
\end{figure}

The number of clusters in a permutation ranges from 1 in the
identity, 2 in cyclic permutations, to $n$ in the majority of the
permutations of $\mathcal S_n$. Let $B(x)$ be the number of clusters
in the permutation of $x$. Observe that the number of clusters
depends only on the element and not on the diagram representing it,
because a reduction or addition of nonreduced carets always happens
inside a cluster, so the number is unchanged. In fact, if the
element of $V$ is represented as a right-continuous map of the
interval $[0,1]$, see \cite{cfp}, the number of clusters is exactly
the number of connected components of the graph of the map. Each
connected component involves the consecutive intervals of a cluster.
See Figure \ref{element}.

With the number of clusters we can state the new, better bound for
the metric for $V$. Let $x\in V$ with $N(x)$ carets and $B(x)$
clusters.

\begin{thm} There exists a constant $K>0$ such that the following
inequality is satisfied for every element $x\in V$:
$$
\frac{N(x)}K\leq ||x||\leq K\left(N(x)+B(x)\log B(x)\right)
$$
\end{thm}

Since $B\leq N$, this bound improves on the bound given by Birget.
Observe that this bound is accurate for elements of $F$ and of $T$.
For $F$, the permutation is the identity, $B(x)=1$, so the bound is
just $N(x)$, as is well-known. For elements of $T$, we have that
$B(x)=2$ and then $B\log B$ is fixed, so the bound is essentially
the number of carets. This bound is, at worst, the same bound as
Birget's, but for a significant number of elements of $V$, namely,
those where $B(x)<N(x)$, it is an improvement.

{\it Proof.} Since the lower bound is the same as Birget's, we only
need to adjust the constant if necessary.

For the upper bound,  let $x$ be an element of $V$, with $N(x)$
carets in each tree of the diagram. We would like to change the
trees, so we pre- and post-multiply it by an element of $F$ (also
with $N(x)$ carets, hence with length $N(x)$ up to the
multiplicative constant $C$), to change the trees arbitrarily. So
without changing significantly the length, we can get an element
whose diagram has the trees we want at both sides of the diagram.
This is the content of the following claim.

\begin{figure}[h!]
\centerline{\includegraphics[width=150mm]{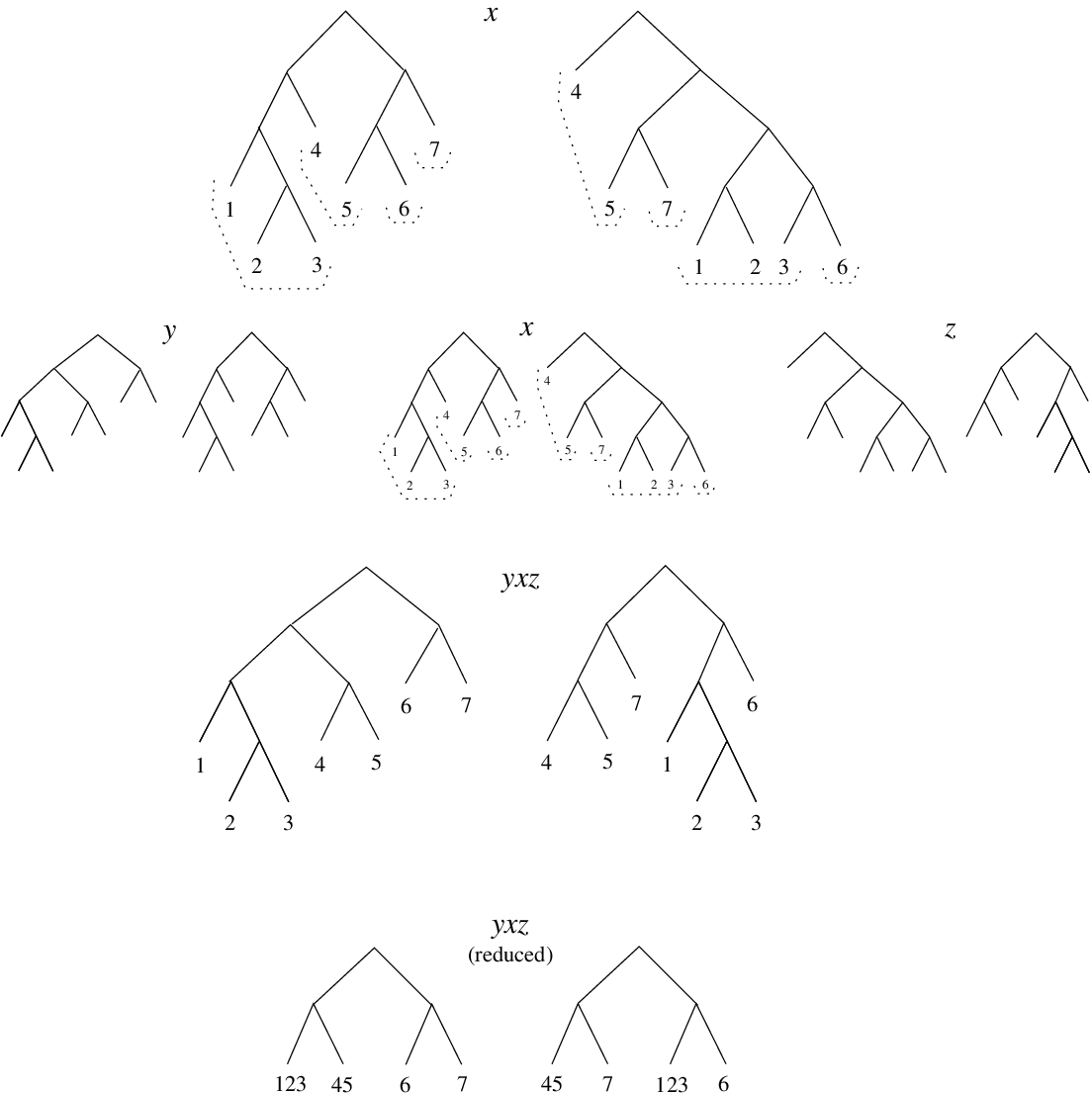}} \caption{The
process of collapsing the clusters of the element in Figure
\ref{element}. In the last diagram each leaf has been labelled with
all the labels of the cluster collapsing into it.} \label{cluster}
\end{figure}

{\it Claim:\/} Given $x\in V$ an element whose reduced diagram  has
$N(x)$ carets and $B(x)$ clusters, there exist two elements $y,z\in
F$, each with a diagram with $N(x)$ carets in their trees, such that
the reduced diagram for $yxz$ has $B(x)$ leaves and hence $B(x)-1$
carets.

Verifying this claim concludes the proof. For the elements $y$ and
$z$ we have $N(y)=N(z)=N(x)$, so they have length approximately
$N(x)$, according to the inequality \ref{eq:2} above. Choose $y$ and
$z$ in such a way that the corresponding diagram for $yxz$ is
nonreduced, and each cluster has the same subtree on both sides,
hence the whole cluster can be reduced to a single leaf. Thus, after
the cancellation, the diagram for $yxz$ has $B(x)-1$ carets and
Birget's bound applies. To find the length for $x$ we combine all
inequalities:

\begin{align} ||x||&\leq ||y||+||yxz||+||z||\\
&\leq C\, N(y) + C\,N(yxz)\log N(yxz) + C\,N(z)\\
&\leq C\, N(x) + C\,(B(x)-1)\log(B(x)-1) + C\,N(x)\\
&\leq 2C\,N(x) + C\,B(x)\log B(x).
\end{align}

and we only need to take $K=2C$ to finish the proof. Observe  that
the constant $K$ also works for the lower bound, so the proof is
complete.  See Figure \ref{cluster} for an example.\qed.

\section {A counterexample to show the bound is not sharp}

\begin{figure}[h!]
\centerline{\includegraphics[width=150mm]{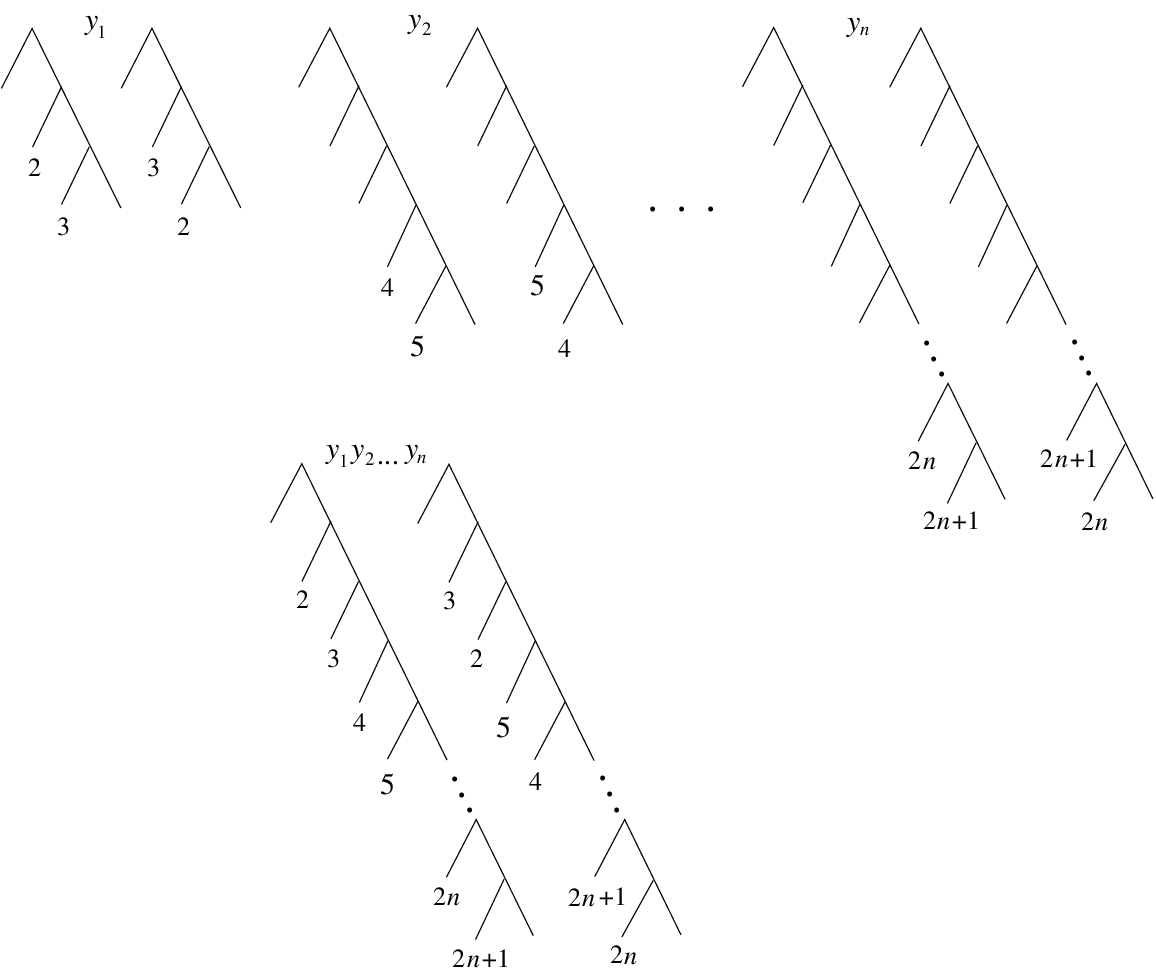}} \caption{The
counterexample to show the bound is not sharp. The unlabelled leaves
are not permuted.} \label{counter}
\end{figure}

Ideally, such a bound would be sharp, meaning that it has the same
quantity on both sides of the inequality (hence providing an
estimate for the distance). Unfortunately, this is not the case,
since $N(x)+B(x)\log B(x)$ is not a lower bound for the metric. It
is not hard to produce an element for $V$ with $N$ carets, and also
with $N$ clusters, but whose distance is also $N$, hence being on
the low side of the bound. Let $y_n$ be the element which has as a
tree just an all right tree with $2n+1$ carets (and hence $2n+2$
leaves), and as permutation the transposition of the leaves $2n$ and
$2n+1$. See Figure \ref{counter} for its diagram.  Clearly
$y_n=x_0^{-2n+2}{y_1}x_0^{2n-2}$. Then, the product
$$
y_1y_2\ldots y_n=y_1x_0^{-2}y_1x_0^{-2}y_1\ldots x_0^{-2}y_1x_0^{2n-2}
$$
has $n$ carets, $n$ clusters, and also length $n$ (all up to a
multiplicative constant).

\bibliography{pepsrefs}

\def\cprime{$'$}
\begin{thebibliography}{1}

\bibitem{birget}
Jean-Camille Birget.
\newblock The groups of {R}ichard {T}hompson and complexity.
\newblock {\em Internat. J. Algebra Comput.}, 14(5-6):569--626, 2004.
\newblock International Conference on Semigroups and Groups in honor of the
  65th birthday of Prof. John Rhodes.

\bibitem{burillo}
Jos{\'e} Burillo.
\newblock Quasi-isometrically embedded subgroups of {T}hompson's group ${F}$.
\newblock {\em J. Algebra}, 212(1):65--78, 1999.

\bibitem{bcs}
Jos\'e Burillo, Sean Cleary, and Melanie Stein.
\newblock Metrics and embeddings of generalizations of {T}hompson's group
  ${F}$.
\newblock {\em Trans. Amer. Math. Soc.}, 353(4):1677--1689 (electronic), 2001.

\bibitem{cfp}
J.~W. Cannon, W.~J. Floyd, and W.~R. Parry.
\newblock Introductory notes on {R}ichard {T}hompson's groups.
\newblock {\em Enseign. Math. (2)}, 42(3-4):215--256, 1996.

\end{thebibliography}
\bibliographystyle{plain}

\end{document}